\documentclass[a4paper]{jpconf}

\providecommand{\keywords}[1]{\textbf{\textit{keywords:}} #1}
\usepackage{float}
\usepackage{amsmath,amsfonts,amssymb,amsthm,epsfig,epstopdf,titling,url,array,tabularx}
\usepackage{caption,subcaption,sublabel}
\usepackage{multirow,multicol}
\usepackage{graphicx}
\usepackage{pgfplots}
\usepackage{algorithmic} 
\usepackage[boxruled,longend,linesnumbered,ruled]{algorithm2e}
\usepackage{cite}
\usepackage{textcomp}
\usepackage{xcolor}
\usepackage{flushend}
\usepackage{eso-pic}
\theoremstyle{definition}

\newtheorem{lemma}{Lemma}[section]

\theoremstyle{remark}

\begin{document}
\title{Interpolatory projection technique for Riccati-based feedback stabilization of index-1 descriptor systems}

\author{Mahtab Uddin$^{1,}$\footnote[5]{Corresponding author}, M Monir Uddin$^2$, M A Hakim Khan$^3$, and Md Motlubar Rahman$^4$}

\address{$^1$Institute of Natural Sciences, United International University, Dhaka-1212, Bangladesh}
\address{$^2$Department of Mathematics \& Physics, North South University, Dhaka-1229, Bangladesh}
\address{$^3$Department of Mathematics, Bangladesh University of Engineering \& Technology, Dhaka-1000, Bangladesh}
\address{$^4$Department of Mathematics, Jahangirnagar University,  Savar, Dhaka-1342, Bangladesh}

\ead{mahtab@ins.uiu.ac.bd, monir.uddin@northsouth.edu, mahkhan@math.buet.ac.bd, mmrmaths@gmail.com}

\begin{abstract}
The work aims to stabilize the unstable index-1 descriptor systems by Riccati-based feedback stabilization via a modified form of Iterative Rational Krylov Algorithm (IRKA), which is a bi-tangential interpolation-based technique. In the basic IRKA, for the stable systems the Reduced Order Models (ROMs) can be found conveniently, but it is unsuitable for the unstable ones. In the proposed technique, the initial feedback is implemented within the construction of the projectors of the IRKA approach. The solution of the Riccati equation is estimated from the ROM achieved by IRKA and hence the low-rank feedback matrix is attained. Using the reverse projecting process, for the full model the optimal feedback matrix is retrieved from the low-rank feedback matrix. Finally, to validate the aptness and competency of the proposed technique it is applied to unstable index-1 descriptor systems. The comparison of the present work with two previous works is narrated. The simulation is done by numerical computation using MATLAB, and both the tabular method and graphical method are used as the supporting tools of comparative analysis.
\\ \\
\end{abstract}

\section{Introduction}\label{sec:introduction}
The index-1 descriptor system of the first-order form can be written with the input-output relations by means of the sparse block-matrices as   
\begin{equation}  \label{eqn:DS_matrix_vector}
\begin{aligned}
\underbrace{ \begin{bmatrix}  E_{1}  &  0 \\ 0  &  0   \end{bmatrix}}_E
\underbrace{ \begin{bmatrix}  \dot{x_1}(t)  \\ \dot{x_2}(t)    \end{bmatrix}}_{\dot{x}(t)}
&= \underbrace{\begin{bmatrix}  J_{1} &  J_{2} \\ J_{3} & J_{4} \end{bmatrix}}_A
\underbrace{ \begin{bmatrix}    x_1(t) \\ x_2(t)  \end{bmatrix}}_{x(t)}
+ \underbrace{\begin{bmatrix}    B_1 \\ B_2  \end{bmatrix}}_B u(t),\\
y(t)&= \underbrace{\begin{bmatrix}  C_1  & C_2 \end{bmatrix}}_C
\begin{bmatrix}    x_1(t) \\ x_2(t)  \end{bmatrix} + D u(t).
\end{aligned}  
\end{equation} 

In system (\ref{eqn:DS_matrix_vector}), $E\in \mathbb{R}^{n\times n}$ and $A\in \mathbb{R}^{n\times n}$ are the differential coefficient matrix and state matrix, respectively, whereas $B\in \mathbb{R}^{n\times p}$ and $C\in \mathbb{R}^{m\times n}$ are the control multiplier matrix and state multiplier matrix, respectively. The matrix $D\in \mathbb{R}^{m\times p}$ is the direct transmission map and it can be absent in some engineering applications, such as power systems models \cite{hossain2019iterative,benner2016reduced}. Here, $n$ is the system dimension, which is very large and $p,m$ are very smaller than $n$. The vectors $x(t)\in \mathbb{R}^n$, $u(t)\in \mathbb{R}^p$, and $y(t)\in \mathbb{R}^m$ represent the state, control (input), and output, respectively, with $x(t_0)=x_0$ as the initial state. The state vector $x(t)$ can be partitioned as $x_1\in \mathbb{R}^{n_1}$, $x_2\in \mathbb{R}^{n_2}$ with $n_1+n_2 = n$. The sub-matrices (blocks) in the system (\ref{eqn:DS_matrix_vector}) are sparse with proper dimensions and $J_4$ is invertible.

For further manipulation, the \texttt{Schur complements} need to be attained from the index-1 descriptor system (\ref{eqn:DS_matrix_vector}) as 
\begin{equation} \label{eqn:schur_complements}
\begin{aligned}
x:&=x_1, \quad 
\mathcal{E}:= E_1, \quad 
\mathcal{A}:= J_1-J_2{J_4}^{-1}J_3, \\  
\mathcal{B}:&= B_1-J_2{J_4}^{-1}B_2, \quad 
\mathcal{C}:= C_1-C_2{J_4}^{-1}J_3, \quad
\mathcal{D}:= D-C_2{J_4}^{-1}B_2.
\end{aligned}
\end{equation} 

Using the block-matrix properties and the \texttt{Schur complements} (\ref{eqn:schur_complements}), index-1 descriptor system (\ref{eqn:DS_matrix_vector}) can be simplified to the generalized LTI continuous-time system as     
\begin{equation} \label{eqn:d-state-space}
\begin{aligned}
\mathcal{E}\dot{x}(t)&=\mathcal{A}x(t)+\mathcal{B}u(t),\\
y(t)&=\mathcal{C}x(t)+\mathcal{D}u(t).
\end{aligned}
\end{equation}

\begin{lemma} [\textbf{Equivalence of transfer functions}] \label{transfer_function}
Assume $G(s)=C(s E - A)^{-1} B + D$  and  $\mathcal{G}(s)=\mathcal{C}(s \mathcal{E} - \mathcal{A})^{-1} \mathcal{B} + \mathcal{D}$ define the transfer functions obtained from the index-1 descriptor system (\ref{eqn:DS_matrix_vector}) and the structured generalized system (\ref{eqn:d-state-space}), respectively. Then, the transfer functions $G(s)$ and $\mathcal{G}(s)$ are identical and because of that those systems are equivalent.
\end{lemma} 

In the topics of applied mathematics and engineering fields, the necessity of LTI continuous-time systems is inevitable, for instance, control theory, system automation, and mechatronics \cite{uddin2015computational,benner2016structure}. The Continuous-time Algebraic Riccati Equation (CARE) emerges within the scientific and engineering applications, for instance, mechanical and electrical systems \cite{chu2011solving,chen2016linear}. The system (\ref{eqn:d-state-space}) generates the CARE formed as follows
\begin{equation} \label{eqn:D-GCARE}
\begin{aligned}
\mathcal{A}^T X\mathcal{E}+\mathcal{E}^T X\mathcal{A}-\mathcal{E}^T X\mathcal{B}\mathcal{B}^T X\mathcal{E}+\mathcal{C}^T \mathcal{C}=0.
\end{aligned}
\end{equation}

The CARE (\ref{eqn:D-GCARE}) is considered as solvable and a unique solution $X$ exists if the Hamiltonian matrix corresponding to the system (\ref{eqn:d-state-space}) has no pure imaginary eigenvalues \cite{abou2012matrix}. For a stable closed-loop matrix $\mathcal{A}- (\mathcal{B}\mathcal{B}^T)X\mathcal{E}$, the symmetric positive-definite $X$ is called stabilizing. 

If the system (\ref{eqn:d-state-space}) is unstable, Riccati-based feedback stabilization is the vital approach and the optimal feedback matrix $K^o$ is the key ingredient, which can be computed as $K^o = \mathcal{B}^T X\mathcal{E}$ \cite{BaeBSetal13}. Using the optimal feedback matrix $K^o$, the target system can be optimally stabilized by replacing $\mathcal{A}$ by $\mathcal{A}_s = \mathcal{A}-\mathcal{B}K^o$. The stabilized system can be written as
\begin{equation} \label{eqn:st-state-space}
\begin{aligned}
\mathcal{E}\dot{x}(t) = \mathcal{A}_s x(t)+\mathcal{B}u(t), \\
y(t) = \mathcal{C}x(t)+\mathcal{D}u(t).
\end{aligned}
\end{equation} 

In some recent works, the Riccati-based feedback stabilization technique for index-1 descriptor systems by the Rational Krylov Subspace Method (RKSM) via the Linear Quadratic Regulator (LQR) approach and the Low-Rank Cholesky-Factor integrated Alternative Direction Implicit (LRCF-ADI) based Kleinman-Newton method ware discussed \cite{uddin2019riccati,uddin2019efficient}. In both of the works, to optimally stabilize, some unstable systems are considered as the target system, those arising from the Brazilian Interconnected Power System (BIPS) models. In this work, we propose a modified form of the bi-tangential interpolation-based Iterative Rational Krylov Algorithm (IRKA) technique for the optimal feedback stabilization of those power system models. A comparative analysis of the proposed method and the previous two methods is included as well.  

\section{IRKA for first-order generalized systems}\label{sec:preliminaries}
In this section, the IRKA approach for the generalized system of first-order is discussed. The first-order representation of the generalized LTI continuous-time system can be formed as  
\begin{equation}\label{eq:ltisystem}
\begin{aligned}
E\dot x(t) &= A x(t)+B u(t),\\
y(t) &= C x(t)+ D u(t),
\end{aligned}
\end{equation}
where $E \in \mathbb{R}^{n\times n}$ is non-singular, and $ A \in \mathbb R ^{n\times n}$, $B \in \mathbb R^{n\times p}$, $ C\in\mathbb R ^{m\times n}$ and ${D}\in\mathbb R^{m\times p}$. 

Assume $\left\{\alpha_i\right\}_{i=1}^{r}\subset\mathbb{C}$ and $\left\{\beta_i\right\}_{i=1}^{r}\subset\mathbb{C}$ are the sets of distinct interpolation points with $b_i \in \mathbb{C}^m$ and $c_i \in \mathbb{C}^p$ as the right and left tangential directions, respectively. These assumptions are essential to construct the left and right projector matrices $W$ and $V$, respectively. The projector matrices $W$ and $V$ can be constructed as  
\begin{equation}\label{eq:solution1}
\begin{aligned}
\text{Range}(V) &=\text{span}\{(\alpha_1 E-A)^{-1}Bb_1, \cdots,(\alpha_r E-A)^{-1}Bb_r\},\\
\text{Range}(W) &=\text{span}\{(\beta_1 E^T-A^T)^{-1}C^Tc_1,\cdots,(\beta_r E^T-A^T)^{-1}C^Tc_r\}.
\end{aligned}
\end{equation}

Employing the projector matrices $W$ and $V$, the reduced-order matrices can be found as follows
\begin{equation}\label{eq:matrixrom}
\begin{aligned}
\hat{E}: &= {W}^T  E V,\quad\quad \hat{A}:= W^T A V,\\
\hat{B}:&= W ^T B , \quad \hat{C}:= C V, \quad \hat{D}:=  D.
\end{aligned}
\end{equation}

Using the above reduced-order matrices, desired Reduced Order Model (ROM) of the system (\ref{eq:ltisystem}) can be written as
\begin{equation}\label{eq:rom}
\begin{aligned}
\hat{E}\dot {\hat{x}}(t) &= \hat{A} \hat{x}(t)+\hat{B} u(t),\\ 
\hat{y}(t) &=\hat{C}\hat{x}(t)+ \hat{D} u(t).
\end{aligned}
\end{equation}

The tangential interpolation projection-based method IRKA was provided in \cite{morXuZ11,morBenS11}, where the authors demonstrated the necessity of the improvement of interpolation points and tangential directions to the desired level.  

For the ROM of size $r$ with $r\ll n$ and $i=1,\cdots, r$, the transfer function $\text{G}(s)$ of the system (\ref{eq:ltisystem}) can be gradually approximated by the transfer function $\hat{\text{G}}(s)$ of the system (\ref{eq:rom}). Then, some predefined conditions need to be fulfilled, such as    
\begin{equation}\label{eq:tfmatrix}
\begin{aligned}
\text{G}(\alpha_i) b_i &= \hat{\text{G}}(\alpha_i) b_i,\\
c_i^T \text{G}(\beta_i) &= c_i^T\hat{\text{G}}(\beta_i), \\
c_i^T \text{G}(\alpha_i) b_i &= c_i^T \hat{\text{G}}(\alpha_i) b_i, \quad \text{when} \ \alpha_i=\beta_i.
\end{aligned}
\end{equation}

For the non-negative integer $j$ and the interpolation point $\alpha_i$, the following Hermite bi-tangential interpolation conditions need to be satisfied 
\begin{equation}\label{eq:multimoment}
\begin{aligned}
c_i^T \text{G}^{(j)}(\alpha_i) b_i &= c_i^T \hat{\text{G}}^{(j)}(\alpha_i) b_i, \\
c_i^T C[(\alpha_i E-A)^{-1} E]^j(\alpha_i E-A)^{-1}B b_i &= c_i^T \hat{C}[(\alpha_i\hat{E}-\hat{A})^{-1}\hat{E}]^j
(\alpha_i\hat{E}-\hat{A})^{-1}\hat{B} b_i,
\end{aligned}
\end{equation}
where $\text{G}^{(j)}(.)$ represents the $j$-th derivative of $\text{G}(.)$ and called the $j$-th moment of $\text{G}(.)$. Algorithm~\ref{alg:irka1} outlines the procedures of the IRKA for the first-order systems.
\begin{algorithm}[]
\SetAlgoLined
\SetKwInOut{Input}{Input}
\SetKwInOut{Output}{Output}
\caption{IRKA for first-order generalized systems.}
\label{alg:irka1}
\Input {$E, A, B, C, D$.}
\Output{$\hat{E}, \hat{A}, \hat{B}, \hat{C}$, $\hat{D}:= D$.}
Assume the initial interpolation points $\{\alpha_i\}_{i=1}^r$ and tangential directions $\{b_i\}_{i=1}^r$ \& $\{c_i\}_{i=1}^r$. \\
Construct $V$ and $W$ according to (\ref{eq:solution1}). \\
\While{(not converged)}{%
$i=1$; \\
Evaluate reduced-order matrices as (\ref{eq:matrixrom}).\\
For $i=1,\cdots ,r$, compute $\hat{A}z_i = \lambda_i\hat{E}z_i$ and $y^*_i \hat{A} = \lambda_i y^*_i \hat{E}$. Then, assign $\alpha_i \leftarrow -\lambda_i$, $b^*_i \leftarrow -y^*_i \hat{B}$ and $c^*_i \leftarrow \hat{C}z^*_i$. \\
Update $V$ and $W$. \\
$i=i+1$;}
Determine the finalized reduced-order matrices by repeating Step-5.\\
\end{algorithm}

\section{Modified IRKA for first-order index-1 descriptor systems}\label{sec:modified}
To modify the basic IRKA, at first the projectors $V$ and $W$ are determined by a set of sparse simplified linear systems. The treatment for the unstable systems is introduced and stabilization of the original system is done by the low-rank feedback matrix attained from ROM.

\subsection{Sparsity preservation}
For the converted generalized system (\ref{eqn:d-state-space}), the projector matrices $W$ and $V$ defined in (\ref{eq:solution1}) need to be written as
\begin{equation}\label{eq:solution2}
\begin{aligned}
\text{Range}(V) &=\text{span}\{(\alpha_1 \mathcal{E}-\mathcal{A})^{-1}Bb_1, 
\cdots,(\alpha_r \mathcal{E}-\mathcal{A})^{-1}Bb_r\},\\
\text{Range}(W) &=\text{span}\{(\alpha_1 \mathcal{E}^T-\mathcal{A}^T)^{-1}C^Tc_1,
\cdots,(\alpha_r \mathcal{E}^T-\mathcal{A}^T)^{-1}C^Tc_r\}.
\end{aligned}
\end{equation}

Due to the dense form of the matrix $\mathcal{A}$, the above forms are usually dense and infeasible for the systems very large dimension. To avoid this adversity, in every iteration two simplified linear systems need to be solved as
\begin{equation} \label{eqn:sparse_shifted1}
\begin{aligned}
(\alpha_i \mathcal{E} - \mathcal{A})v_i &=\mathcal{B} b_i,\\
or, {\begin{bmatrix} (\alpha_i E_1 - J_1)  &  -J_2 \\ -J_3  &  -J_4 \end{bmatrix}} \begin{bmatrix}  v_i  \\ \Gamma_1 \end{bmatrix}
&=  \begin{bmatrix}    B_1 \\ B_2  \end{bmatrix}b_i,
\end{aligned}
\end{equation}
and
\begin{equation} \label{eqn:sparse_shifted2}
\begin{aligned}
(\alpha_i \mathcal{E}^T - \mathcal{A}^T)w_i &=\mathcal{C}^T c_i,\\
or, {\begin{bmatrix} (\alpha_i E_1 - J_1)  &  -J_2 \\ -J_3  &  -J_4   \end{bmatrix}}^T \begin{bmatrix}  w_i  \\ \Gamma_2 \end{bmatrix}
&=  \begin{bmatrix}    C_1^T \\ C_2^T  \end{bmatrix}c_i.
\end{aligned}
\end{equation} 

Here, $\Gamma_1$ and $\Gamma_2$ are the truncated terms. The explicit form of the reduced-order matrices in (\ref{eq:rom}) defined in (\ref{eq:matrixrom}) are inefficient to construct feasible ROM. The sparsity preserving reduced-order matrices can be attained by following way
\begin{equation} \label{eqn:matrixsprom}
\begin{aligned}
\hat{\mathcal{E}} &= W^T E_1 V,\quad
\hat{\mathcal{A}} = W^T J_1 V - (W^T J_2){J_4}^{-1}(J_3 V),\\
\hat{\mathcal{B}} &= W^T B_1 - (W^T J_2) {J_4}^{-1} B_2,\quad
\hat{\mathcal{C}} = C_1 V - C_2 {J_4}^{-1} (J_3 V).
\end{aligned}
\end{equation} 

\subsection{Treatment for the unstable systems}
In the case of unstable index-1 descriptor systems, the Bernoulli stabilization process is required. Using an initial feedback matrix $K_0$, the matrix $\mathcal{A}$ needs to be replaced by $\mathcal{A}_f = \mathcal{A}-\mathcal{B} K_0$. Then, the re-structured form of system (\ref{eqn:d-state-space}) and corresponding CARE (\ref{eqn:D-GCARE}) can be written as
\begin{equation} \label{eqn:re-state-space}
\begin{aligned}
\mathcal{E}\dot{x}(t)&=\mathcal{A}_f x(t)+\mathcal{B}u(t),\\
y(t)&=\mathcal{C}x(t)+\mathcal{D}u(t),
\end{aligned}
\end{equation}
\begin{equation} \label{eqn:re-GCARE}
\begin{aligned}
\mathcal{A}_f^T X\mathcal{E}+\mathcal{E}^T X\mathcal{A}_f-\mathcal{E}^T X\mathcal{B} \mathcal{B}^T X\mathcal{E}+\mathcal{C}^T \mathcal{C}=0.
\end{aligned}
\end{equation}  

Then, the simplified linear systems (\ref{eqn:sparse_shifted1}) and (\ref{eqn:sparse_shifted2}) can be re-written as 
\begin{equation} \label{eqn:sparse_shifted3}
\begin{aligned}
(\alpha_i \mathcal{E} - \mathcal{A}_f)v_i &=\mathcal{B}b_i,\\
or, {\begin{bmatrix} (\alpha_i E_1 - (J_1 - B_1 K_0))  &  -J_2 \\ -(J_3 - B_2 K_0)  &  -J_4   \end{bmatrix}} \begin{bmatrix}  v_i  \\ \Gamma_1 \end{bmatrix}
&=  \begin{bmatrix}    B_1 \\ B_2  \end{bmatrix}b_i,
\end{aligned}
\end{equation}
and 
\begin{equation} \label{eqn:sparse_shifted4}
\begin{aligned}
(\alpha_i \mathcal{E}^T - \mathcal{A}_f^T)w_i &=\mathcal{C}^T c_i,\\
or, {\begin{bmatrix} (\alpha_i E_1 - (J_1 - B_1 K_0))  &  -J_2 \\ -(J_3 - B_2 K_0)  &  -J_4   \end{bmatrix}}^T \begin{bmatrix}  w_i  \\ \Gamma_2 \end{bmatrix}
&=  \begin{bmatrix}    C_1^T \\ C_2^T  \end{bmatrix}c_i.
\end{aligned}  
\end{equation} 

\subsection{Optimal feedback matrix from the ROMs}
The ultimate goal of the work is included here. Using the reduced-order matrices defined in (\ref{eqn:matrixsprom}), the reduced-order form of the system (\ref{eqn:d-state-space}) and corresponding CARE can be found as
\begin{equation} \label{eqn:r-state-space}
\begin{aligned}
\hat{\mathcal{E}}\dot{x}(t)&=\hat{\mathcal{A}}x(t)+\hat{\mathcal{B}}u(t),\\
y(t)&=\hat{\mathcal{C}}x(t)+\hat{\mathcal{D}}u(t),
\end{aligned}
\end{equation}
\begin{equation} \label{eqn:lr-GCARE}
\begin{aligned}
\hat{\mathcal{A}}^T \hat{X}\hat{\mathcal{E}} + \hat{\mathcal{E}}\hat{X}\hat{\mathcal{A}} - \hat{\mathcal{E}}\hat{X}\hat{\mathcal{B}}\hat{\mathcal{B}}^T \hat{X}\hat{\mathcal{E}} -\hat{\mathcal{C}}^T\hat{\mathcal{C}} =0.
\end{aligned}
\end{equation}

Now, the optimal feedback matrix for the original model can be approximated by employing the ROM (\ref{eqn:r-state-space}). For this, the low-rank CARE (\ref{eqn:lr-GCARE}) needs to be solved for symmetric positive-definite matrix $\hat{X}$ by the MATLAB library command \texttt{care}. Then, the stabilizing feedback matrix $\hat{K} = \hat{\mathcal{B}}^T\hat{X}$ for the ROM (\ref{eqn:r-state-space}) can be estimated, and hence for stabilizing the full model (\ref{eqn:d-state-space}) the approximated optimal feedback matrix $K^o$ can be retrieved as
\begin{equation}
\begin{aligned}
K^o = \hat{\mathcal{B}}^T\hat{X}{V}^TE_1 = \hat{K}{V}^TE_1,
\end{aligned}
\end{equation}
where $V$ is the right projector matrix. Algorithm~\ref{alg:irka2} outlines the modified form of the IRKA for the first-order index-1 descriptor systems.  
\begin{algorithm}[]
\SetAlgoLined
\SetKwInOut{Input}{Input}
\SetKwInOut{Output}{Output}
\caption{Modified IRKA for first-order index-1 descriptor systems.}
\label{alg:irka2}
\Input {$E_1, J_1, J_2, J_3, J_4, B_1, B_2, C_1, C_2, D$.}
\Output{$\hat{\mathcal{E}}, \hat{\mathcal{A}}, \hat{\mathcal{B}}, \hat{\mathcal{C}}, \hat{D}:= D$ and optimal feedback matrix $K^o$.}
Assume the initial interpolation points $\{\alpha_i\}_{i=1}^r$ and tangential directions $\{b_i\}_{i=1}^r$ \& $\{c_i\}_{i=1}^r$. \\
Construct $V = \begin{bmatrix}v_1, v_2, \cdots, v_r \end{bmatrix}$ and
$W = \begin{bmatrix}w_1, w_2, \cdots, w_r \end{bmatrix}$, where $v_i$ and $w_i$ can be attained from (\ref{eqn:sparse_shifted3}) and (\ref{eqn:sparse_shifted4}). \\
\While{(not converged)}{%
$i=1$; \\
Evaluate sparsity preserving $\hat{\mathcal{E}}$, $\hat{\mathcal{A}}$, $\hat{\mathcal{B}}$, and $\hat{\mathcal{C}}$ by (\ref{eqn:matrixsprom}).\\ 
For $i=1,\cdots ,r$, compute $\hat{\mathcal{A}}z_i = \lambda_i \hat{\mathcal{E}}z_i$ and $y^*_i \hat{\mathcal{A}} = \lambda_i y^*_i \hat{\mathcal{E}}$. Then, assign $\alpha_i \leftarrow -\lambda_i$, $b^*_i \leftarrow -y^*_i \hat{\mathcal{B}}$ and $c^*_i \leftarrow \hat{\mathcal{C}}z^*_i$. \\
Update $V$ and $W$. \\
$i=i+1$;}
Determine the finalized reduced-order matrices by repeating Step-5.\\
Solve the low-rank Riccati equation (\ref{eqn:lr-GCARE}) for $\hat{X}$.\\
Compute the reduced-order feedback matrix $\hat{K} = \hat{\mathcal{B}}^T\hat{X}$.\\
Retrieve the optimal feedback matrix $K^o = \hat{K}{V}^TE_1$.
\end{algorithm}

At the final step of the computation, using the optimal feedback matrix $K^o$, the system (\ref{eqn:d-state-space}) can be optimally stabilized as (\ref{eqn:st-state-space}).

\section{Numerical results}\label{sec:num_results}
\newlength\figwidth
\setlength{\figwidth}{.33\linewidth}
\newlength\figheight
\setlength{\figheight}{.2\linewidth}
\tikzset{mark options={solid,mark size=3,line width=5pt,mark repeat=10},line width=5pt}

To justify the validity of the proposed technique, it is applied to the data generated in some large-scale real-world models. The computations are carried out with MATLAB\textsuperscript{\textregistered} R2015a (8.5.0.197613) on a board with Intel\textsuperscript{\textregistered}$ \text{Core}^{\text{TM}}$i5 6200U CPU with a 2.30~GHz clock speed and 16~GB RAM. For numerical computations, the following model examples are used.

\subsection{Brazilian Interconnected Power System (BIPS) models}
Brazilian Interconnected Power System (BIPS) models are the most convenient examples of the index-1 descriptor systems \cite{freitas1999computationally}. The following Table.\ref{tab:BIPS_systems} provides the details about the unstable power system models\footnote[6]{https://sites.google.com/site/rommes/software} $mod-606$, $mod-1998$, and $mod-2476$, where the names of the models are considered according to their number of states \cite{leandro2015identification}.
\begin{table}[]
\centering
\caption{Structure of the unstable power system models}
\label{tab:BIPS_systems}
\begin{tabular}{cccc}\hline
Dimensions &  7135   &   15066    &   16861 \\ \hline
States &  606   &   1998    &   2476 \\ \hline
Algebraic variables &   6529  &   13068    &   14385 \\ \hline
Inputs &  4   &    4   &    4   \\ \hline
Outputs &  4   &   4    &    4   \\ \hline
\end{tabular}
\end{table}

\subsection{Reduced-order models of BIPS models}
To find the stabilizing feedback matrix for unstable BIPS model, at first the ROMs need to be computed. The ROMs of the models $mod-606$, $mod-1998$, and $mod-2476$ are computed with the dimensions $30$, $70$, and $100$, respectively. In the numerical computations, for every model the truncation tolerance $10^{-5}$ for the relative error and maximum number of iterations $i_{max}=150$ have been taken.  For the convenience of time and compactness of this work, the analysis of the results found for the model $mod-2476$ will be narrated.

\begin{figure}[]
\setlength{\figheight}{.325\linewidth}
\begin{subfigure}{\linewidth}
\centering
\setlength{\figwidth}{.75\linewidth}
\input{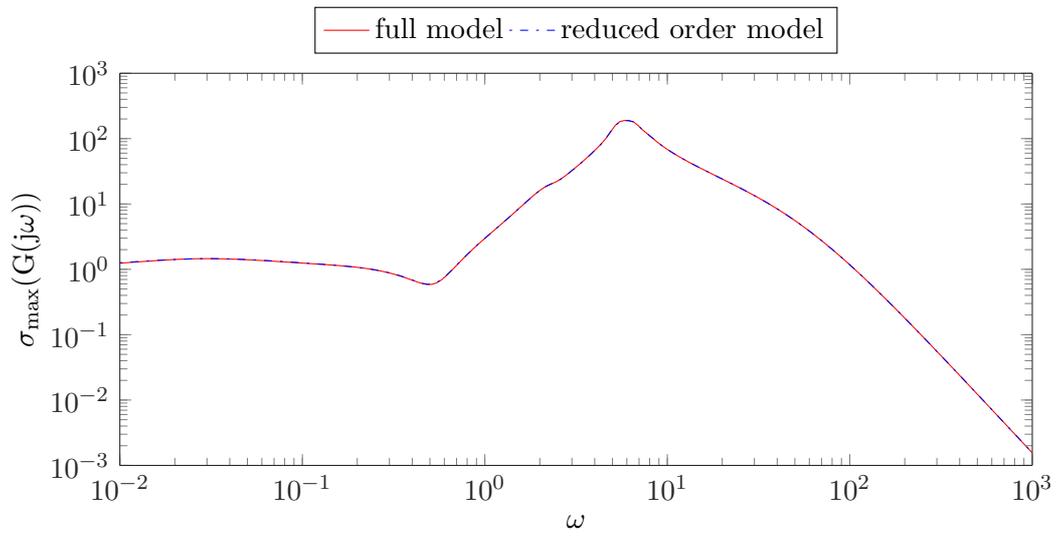}
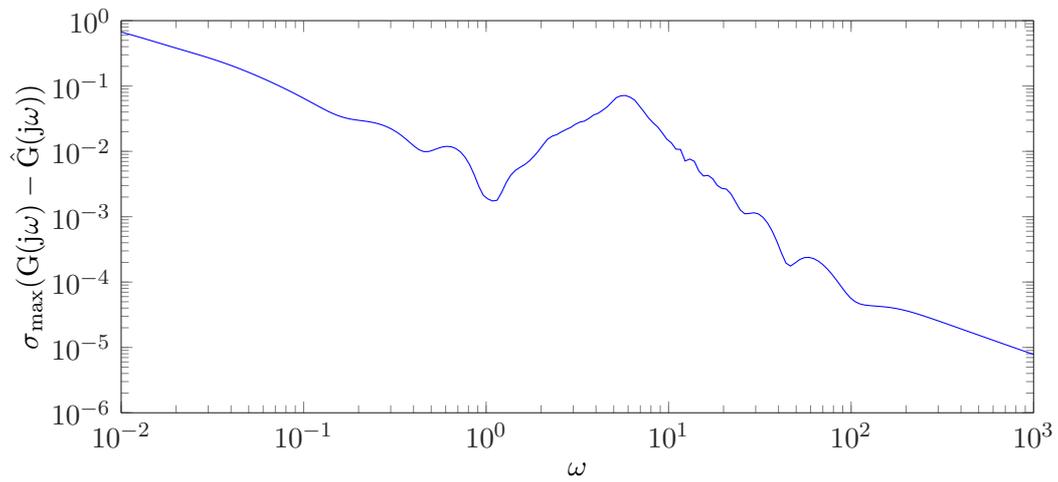
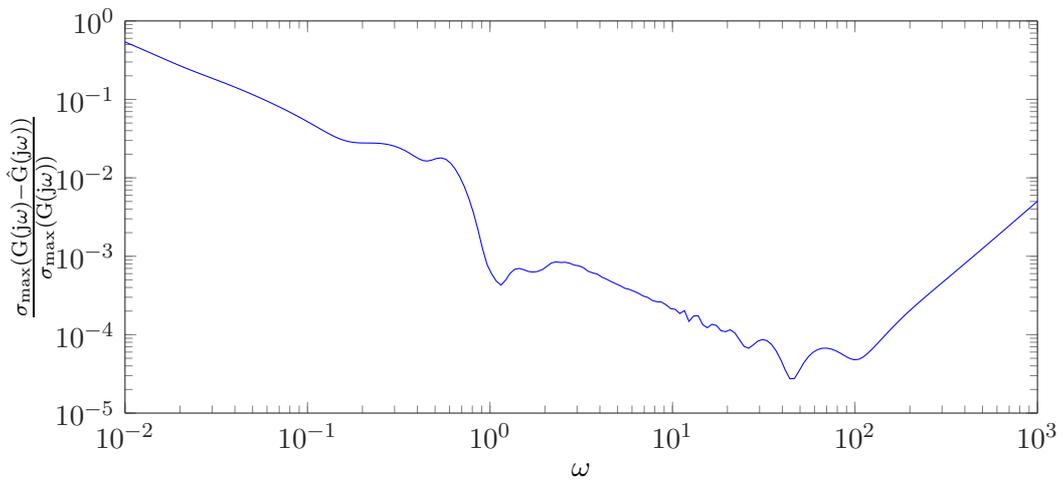
\caption{Sigma plot}
\label{fig:sigmaplot}
\end{subfigure}
\begin{subfigure}{\linewidth}
\centering
\setlength{\figwidth}{.75\linewidth}
%
%
%
%
\begin{tikzpicture}

\begin{axis}[%
width=\figwidth,
height=\figheight,
scale only axis,
separate axis lines,
every outer x axis line/.append style={darkgray!60!black},
every x tick label/.append style={font=\color{darkgray!60!black}},
xmode=log,
xmin=0.01,
xmax=1000,
xminorticks=true,
xlabel={$\omega$},
every outer y axis line/.append style={darkgray!60!black},
every y tick label/.append style={font=\color{darkgray!60!black}},
ymode=log,
ymin=1e-06,
ymax=1,
yminorticks=true,
ylabel={${\sigma{}_{\text{max}}\text{(G(j}\omega\text{)}-\hat{\text{G}}\text{(j}\omega\text{))}}$},
]
\addplot [
color=blue,
solid,
]
table[row sep=crcr]{
0.01 0.669855537797067\\
0.0105956017927762 0.640653723015392\\
0.0112266777351081 0.612214189387393\\
0.0118953406737032 0.584592923129048\\
0.0126038292967973 0.557842281081406\\
0.013354515629299 0.532011895894069\\
0.0141499129743458 0.507148225894957\\
0.0149926843278605 0.483292664945812\\
0.0158856512942805 0.460478337614642\\
0.0168318035333096 0.438725942941625\\
0.0178343087693191 0.418039339678148\\
0.0188965233969121 0.398401747716188\\
0.0200220037181558 0.379773542059011\\
0.0212145178491063 0.362092436011851\\
0.0224780583354873 0.34527646338705\\
0.0238168555197616 0.329229586756239\\
0.0252353917043477 0.313849206095847\\
0.0267384161583995 0.299034405615451\\
0.0283309610183932 0.284693685223534\\
0.0300183581357559 0.270751088946972\\
0.0318062569279412 0.257150057923931\\
0.0337006432927193 0.24385481442614\\
0.0357078596490046 0.230849495010005\\
0.0378346261713193 0.218135539048481\\
0.0400880632889846 0.205727943603277\\
0.042475715525369 0.193650995044678\\
0.045005576757005 0.181933985395023\\
0.0476861169771447 0.170607294935909\\
0.0505263106533568 0.159699097337262\\
0.0535356667741072 0.149232831165957\\
0.0567242606849198 0.139225500804183\\
0.0601027678207038 0.129686803057608\\
0.0636824994471859 0.12061903155436\\
0.0674754405311069 0.112017674731212\\
0.0714942898659758 0.103872596915324\\
0.0757525025877191 0.096169674205294\\
0.0802643352225717 0.0888927476644704\\
0.0850448934180268 0.0820257534185184\\
0.0901101825166502 0.0755548919219889\\
0.0954771611420806 0.0694707016639131\\
0.101163797976621 0.0637699056164567\\
0.107189131920513 0.0584568872737675\\
0.113573335834311 0.0535445893440759\\
0.120337784077759 0.0490544761929219\\
0.127505124071301 0.0450149921245174\\
0.135099352119803 0.0414577842038791\\
0.143145893752348 0.0384110458035283\\
0.151671688847092 0.0358900037393979\\
0.160705281826164 0.0338860594185059\\
0.17027691722259 0.0323579603988604\\
0.180418640939207 0.031229142045134\\
0.19116440753857 0.0303935637155969\\
0.202550193923067 0.0297285839328146\\
0.21461411978584 0.029110470315228\\
0.227396575235793 0.0284280615376083\\
0.240940356023953 0.0275921664602036\\
0.255290806823952 0.0265404825529367\\
0.270495973046314 0.0252390234183601\\
0.286606761694825 0.0236813686876709\\
0.303677111803546 0.0218870020875812\\
0.321764175025074 0.0198999409706875\\
0.340928506974681 0.0177888954397157\\
0.361234269970943 0.0156502169048984\\
0.382749447851631 0.0136142225925985\\
0.405546073584083 0.0118514983654627\\
0.429700470432084 0.0105624707093534\\
0.455293507486695 0.00991453878675084\\
0.482410870416537 0.00992562284691383\\
0.511143348344017 0.010412209066506\\
0.541587137807947 0.0110834233317624\\
0.57384416483024 0.0116614396025997\\
0.608022426164942 0.0119320035940393\\
0.644236350872137 0.0117507596001353\\
0.682607183427239 0.0110440896361393\\
0.723263389648354 0.00981647111085431\\
0.766341086800746 0.00816110467657851\\
0.811984499318401 0.00626550321679145\\
0.860346441668451 0.00440566019309023\\
0.911588829975082 0.00292475662623205\\
0.96588322411587 0.00212047581205199\\
1.02341140210545 0.00186802996634938\\
1.08436596868961 0.0017404451762469\\
1.14895100018731 0.00178624613269783\\
1.21738272773966 0.00236866951947917\\
1.28989026125331 0.00334960064262684\\
1.36671635646201 0.00434610328263157\\
1.44811822767453 0.00513075405658244\\
1.53436840893001 0.00570163306633659\\
1.62575566644379 0.00627360760467204\\
1.72258596539879 0.00712544435355531\\
1.82518349431904 0.00838179347598402\\
1.93389175045523 0.0101129506321988\\
2.04907468981585 0.0125545340639738\\
2.1711179456945 0.0153570067763591\\
2.30043011977292 0.0170864194293019\\
2.43744415012222 0.0180252129942832\\
2.58261876068267 0.0198417369749645\\
2.73643999707467 0.0215382479840026\\
2.89942285388288 0.0231592633943315\\
3.07211299886176 0.0258459820451525\\
3.25508859983506 0.0279666428967241\\
3.44896226040576 0.0289933850829686\\
3.65438307095726 0.0318473754752726\\
3.87203878181256 0.0358945309888294\\
4.10265810582719 0.0382277561436041\\
4.34701315812503 0.0425412917026069\\
4.60592204114511 0.0480006168800537\\
4.88025158365443 0.0564738693514701\\
5.17092024289676 0.0664481202273416\\
5.47890117959394 0.070919389771031\\
5.8052255160949 0.0715657426867864\\
6.1509857885805 0.0673579980302902\\
6.51733960488243 0.060954757903581\\
6.90551352016233 0.0497110218929627\\
7.3168071434272 0.0406318490671882\\
7.75259748862946 0.0326365611103493\\
8.21434358491943 0.0274455931385052\\
8.70359136148517 0.0237073451999118\\
9.22197882333433 0.0192807566140753\\
9.7712415353465 0.0153853741637665\\
10.3532184329566 0.0135162596383333\\
10.9698579789238 0.0108123320127012\\
11.6232246867985 0.0107058218007927\\
12.3155060329283 0.00711866019142339\\
13.049019780144 0.00761986496992613\\
13.8262217376466 0.00706199412131731\\
14.6497139830729 0.00499528735741889\\
15.5222535742705 0.00420425706062215\\
16.4467617799466 0.00429173636250519\\
17.4263338600965 0.00382211546336048\\
18.4642494289554 0.00301344878354697\\
19.5639834351706 0.00271624765701746\\
20.7292177959537 0.00264666436319875\\
21.9638537241655 0.00222940438386533\\
23.2720247896041 0.00168157360267557\\
24.658110758226 0.00127373385466994\\
26.1267522556333 0.0011090907230127\\
27.6828663039207 0.00112180204139089\\
29.3316627839005 0.00114864014831208\\
31.0786618778201 0.00110314325749987\\
32.9297125509715 0.000976443374758729\\
34.8910121340677 0.000795896297626116\\
36.9691270719503 0.000599894806078813\\
39.1710149080926 0.000420267648247948\\
41.5040475785048 0.000279055660630845\\
43.9760360930272 0.000194670384324878\\
46.5952566866468 0.000176100753902038\\
49.37047852839 0.000196525895325174\\
52.3109930805626 0.000221503649708412\\
55.4266452066311 0.000236596755980993\\
58.7278661318948 0.000238765709526223\\
62.2257083673023 0.000229170345830779\\
65.9318827133355 0.000210563348164211\\
69.8587974678525 0.000186205889082318\\
74.0195999691564 0.000159235376923099\\
78.4282206133768 0.000132313280943098\\
83.099419493534 0.00010750268344757\\
88.0488358164346 8.62842604733049e-05\\
93.2930402628469 6.95986934885497e-05\\
98.8495904662559 5.78012010600586e-05\\
104.737089795945 5.05110772606889e-05\\
110.975249641207 4.66226203896906e-05\\
117.584955405216 4.47550458005119e-05\\
124.588336429501 4.3816708511127e-05\\
132.008840083142 4.32028957061015e-05\\
139.871310264724 4.26597896737134e-05\\
148.202070579886 4.20806922525258e-05\\
157.029012472938 4.13904080766894e-05\\
166.381688607613 4.05270490414711e-05\\
176.291411809595 3.9459821156501e-05\\
186.791359902078 3.81959562835444e-05\\
197.916686785356 3.67706479115897e-05\\
209.704640132323 3.52315984099808e-05\\
222.194686093952 3.36268605913559e-05\\
235.428641432242 3.19981124890646e-05\\
249.450813523032 3.03782001221478e-05\\
264.30814869741 2.87911817539115e-05\\
280.050389418363 2.72535275208348e-05\\
296.730240818887 2.57756712605705e-05\\
314.40354715915 2.43635083908857e-05\\
333.129478793467 2.30196712879917e-05\\
352.970730273065 2.1744539512779e-05\\
373.993730247879 2.05370016229453e-05\\
396.268863870148 1.93950082891544e-05\\
419.870708444391 1.83159610854568e-05\\
444.878283112759 1.72969775733947e-05\\
471.375313411672 1.63350664634538e-05\\
499.450511585514 1.54272393095509e-05\\
529.197873595844 1.45705786156162e-05\\
560.716993820546 1.37622767731083e-05\\
594.113398496503 1.29996560127554e-05\\
629.498899022189 1.22801763619794e-05\\
666.991966303012 1.16014362956522e-05\\
706.718127392749 1.09611691445707e-05\\
748.810385759002 1.03572372131598e-05\\
793.409666579749 9.7876248138847e-06\\
840.665288561832 9.25043093683985e-06\\
890.735463861044 8.743861967417e-06\\
943.787827777537 8.2662246660163e-06\\
1000 7.81591951504204e-06\\
};
\end{axis}
\end{tikzpicture}%
\caption{Absolute error}
\label{fig:abserror}
\end{subfigure}
\begin{subfigure}{\linewidth}
\centering
\setlength{\figwidth}{.75\linewidth}
%
%
%
%
\begin{tikzpicture}

\begin{axis}[%
width=\figwidth,
height=\figheight,
scale only axis,
separate axis lines,
every outer x axis line/.append style={darkgray!60!black},
every x tick label/.append style={font=\color{darkgray!60!black}},
xmode=log,
xmin=0.01,
xmax=1000,
xminorticks=true,
xlabel={$\omega$},
every outer y axis line/.append style={darkgray!60!black},
every y tick label/.append style={font=\color{darkgray!60!black}},
ymode=log,
ymin=1e-05,
ymax=1,
yminorticks=true,
ylabel={$\frac{\sigma{}_{\text{max}}\text{(G(j}\omega\text{)}-\hat{\text{G}}\text{(j}\omega\text{))}}{\sigma{}_{\text{max}}\text{(G(j}\omega\text{))}}$},
]
\addplot [
color=blue,
solid,
]
table[row sep=crcr]{
0.01 0.539103902030671\\
0.0105956017927762 0.509745172670643\\
0.0112266777351081 0.481572455219803\\
0.0118953406737032 0.454598785219277\\
0.0126038292967973 0.428842162879237\\
0.013354515629299 0.404324416018834\\
0.0141499129743458 0.38106887326756\\
0.0149926843278605 0.359096930747124\\
0.0158856512942805 0.338423778907401\\
0.0168318035333096 0.319053738131284\\
0.0178343087693191 0.300975850114731\\
0.0188965233969121 0.284160435865002\\
0.0200220037181558 0.268557301256604\\
0.0212145178491063 0.254096062085344\\
0.0224780583354873 0.240688719787749\\
0.0238168555197616 0.228234211351007\\
0.0252353917043477 0.216624305602866\\
0.0267384161583995 0.205749983644376\\
0.0283309610183932 0.195507438077347\\
0.0300183581357559 0.185802971601659\\
0.0318062569279412 0.176556361930772\\
0.0337006432927193 0.167702563437779\\
0.0357078596490046 0.159191868413421\\
0.0378346261713193 0.150988826352392\\
0.0400880632889846 0.143070286593642\\
0.042475715525369 0.135422937006799\\
0.045005576757005 0.128040660135983\\
0.0476861169771447 0.120921961665632\\
0.0505263106533568 0.114067658436142\\
0.0535356667741072 0.107478950340411\\
0.0567242606849198 0.101155953962378\\
0.0601027678207038 0.0950967350301771\\
0.0636824994471859 0.089296847844827\\
0.0674754405311069 0.0837493633463125\\
0.0714942898659758 0.0784453450965479\\
0.0757525025877191 0.0733747134369775\\
0.0802643352225717 0.0685274229592846\\
0.0850448934180268 0.063894866159616\\
0.0901101825166502 0.059471407086851\\
0.0954771611420806 0.0552559400548544\\
0.101163797976621 0.0512533613373799\\
0.107189131920513 0.0474758263895582\\
0.113573335834311 0.0439436102169153\\
0.120337784077759 0.0406852588878512\\
0.127505124071301 0.0377365325452946\\
0.135099352119803 0.0351374726834599\\
0.143145893752348 0.0329269478866655\\
0.151671688847092 0.0311345313030372\\
0.160705281826164 0.0297707899579372\\
0.17027691722259 0.0288187195397694\\
0.180418640939207 0.0282299122601529\\
0.19116440753857 0.0279277434448768\\
0.202550193923067 0.0278167385972008\\
0.21461411978584 0.027794647032297\\
0.227396575235793 0.0277634735435489\\
0.240940356023953 0.0276373585631881\\
0.255290806823952 0.0273470620354719\\
0.270495973046314 0.0268418760489933\\
0.286606761694825 0.026090114116464\\
0.303677111803546 0.0250793754146308\\
0.321764175025074 0.0238179660241564\\
0.340928506974681 0.0223394136470226\\
0.361234269970943 0.020712952558412\\
0.382749447851631 0.0190632245512084\\
0.405546073584083 0.0175971078534087\\
0.429700470432084 0.0166089553985038\\
0.455293507486695 0.0163789738963313\\
0.482410870416537 0.0169045975561071\\
0.511143348344017 0.0176898026308463\\
0.541587137807947 0.0179606970749901\\
0.57384416483024 0.0171943975116085\\
0.608022426164942 0.0154098373473618\\
0.644236350872137 0.012982970663953\\
0.682607183427239 0.0103300184359281\\
0.723263389648354 0.00775977451793691\\
0.766341086800746 0.00546956826572273\\
0.811984499318401 0.00357942667934135\\
0.860346441668451 0.00215814735317284\\
0.911588829975082 0.00123509233690676\\
0.96588322411587 0.000775222556176426\\
1.02341140210545 0.000592946844166394\\
1.08436596868961 0.000480400159510041\\
1.14895100018731 0.000428907460425141\\
1.21738272773966 0.000494529943280386\\
1.28989026125331 0.000607413443331258\\
1.36671635646201 0.000683611699242203\\
1.44811822767453 0.000699057004605905\\
1.53436840893001 0.000672126122935762\\
1.62575566644379 0.00063945901891499\\
1.72258596539879 0.000628248005218773\\
1.82518349431904 0.000640890127140606\\
1.93389175045523 0.000675110209764482\\
2.04907468981585 0.0007420036810068\\
2.1711179456945 0.000821829312349711\\
2.30043011977292 0.000847265857431175\\
2.43744415012222 0.000832745561420321\\
2.58261876068267 0.000838788557022232\\
2.73643999707467 0.000813703562440473\\
2.89942285388288 0.000771015335446523\\
3.07211299886176 0.000754370424073768\\
3.25508859983506 0.00071366911249049\\
3.44896226040576 0.000644105553747679\\
3.65438307095726 0.000612270609752038\\
3.87203878181256 0.000593600801368182\\
4.10265810582719 0.000541041593137525\\
4.34701315812503 0.000510988353601949\\
4.60592204114511 0.000476517321326477\\
4.88025158365443 0.000446595510197913\\
5.17092024289676 0.000421831895175972\\
5.47890117959394 0.000390629991538376\\
5.8052255160949 0.000378113832717649\\
6.1509857885805 0.000357302930424438\\
6.51733960488243 0.000337077056573234\\
6.90551352016233 0.000311315800836021\\
7.3168071434272 0.000298746242015334\\
7.75259748862946 0.000273502568564807\\
8.21434358491943 0.000262973194793312\\
8.70359136148517 0.000260582683221941\\
9.22197882333433 0.000240409080314023\\
9.7712415353465 0.000215007892103227\\
10.3532184329566 0.000210154083465907\\
10.9698579789238 0.000185931946783761\\
11.6232246867985 0.000202729486388165\\
12.3155060329283 0.000147862507636389\\
13.049019780144 0.000172906516034677\\
13.8262217376466 0.00017454216493098\\
14.6497139830729 0.000134193549674301\\
15.5222535742705 0.00012259019353857\\
16.4467617799466 0.000135713287208092\\
17.4263338600965 0.000131006360041829\\
18.4642494289554 0.000111926227634905\\
19.5639834351706 0.000109320132301473\\
20.7292177959537 0.000115447045526951\\
21.9638537241655 0.00010544016155059\\
23.2720247896041 8.62851081441753e-05\\
24.658110758226 7.09703232606316e-05\\
26.1267522556333 6.71789815584793e-05\\
27.6828663039207 7.39707767057379e-05\\
29.3316627839005 8.25911643257289e-05\\
31.0786618778201 8.66618468855705e-05\\
32.9297125509715 8.39911771670728e-05\\
34.8910121340677 7.51397882230096e-05\\
36.9691270719503 6.23206776420113e-05\\
39.1710149080926 4.81738723915303e-05\\
41.5040475785048 3.53952857600099e-05\\
43.9760360930272 2.74035409751381e-05\\
46.5952566866468 2.75952393034636e-05\\
49.37047852839 3.43868692868646e-05\\
52.3109930805626 4.34107760915942e-05\\
55.4266452066311 5.20972353302791e-05\\
58.7278661318948 5.92524325663758e-05\\
62.2257083673023 6.42904796135731e-05\\
65.9318827133355 6.69777887856603e-05\\
69.8587974678525 6.73567691589928e-05\\
74.0195999691564 6.56925800245301e-05\\
78.4282206133768 6.24284833574722e-05\\
83.099419493534 5.81668420649587e-05\\
88.0488358164346 5.36779608260839e-05\\
93.2930402628469 4.99067672777816e-05\\
98.8495904662559 4.78880006882517e-05\\
104.737089795945 4.84613555144122e-05\\
110.975249641207 5.19111853286029e-05\\
117.584955405216 5.79489580847954e-05\\
124.588336429501 6.61019150280989e-05\\
132.008840083142 7.60740149902187e-05\\
139.871310264724 8.78248720106365e-05\\
148.202070579886 0.000101445257328974\\
157.029012472938 0.000117010129582659\\
166.381688607613 0.000134530417682737\\
176.291411809595 0.000153997442243955\\
186.791359902078 0.000175447169414077\\
197.916686785356 0.000198997437161926\\
209.704640132323 0.000224855633162275\\
222.194686093952 0.000253311427512145\\
235.428641432242 0.000284726518586527\\
249.450813523032 0.000319526904538111\\
264.30814869741 0.000358199103821458\\
280.050389418363 0.000401290084569907\\
296.730240818887 0.000449410262818035\\
314.40354715915 0.000503238995724054\\
333.129478793467 0.000563532171063727\\
352.970730273065 0.000631131657957025\\
373.993730247879 0.000706976514238061\\
396.268863870148 0.000792115944781585\\
419.870708444391 0.000887724084064907\\
444.878283112759 0.000995116742474015\\
471.375313411672 0.00111577031644996\\
499.450511585514 0.00125134311990162\\
529.197873595844 0.00140369945370095\\
560.716993820546 0.00157493679048127\\
594.113398496503 0.00176741651881941\\
629.498899022189 0.00198379876272873\\
666.991966303012 0.00222708187298374\\
706.718127392749 0.00250064727725585\\
748.810385759002 0.00280831047841684\\
793.409666579749 0.00315437910672721\\
840.665288561832 0.00354371906269526\\
890.735463861044 0.00398182993750896\\
943.787827777537 0.00447493106336941\\
1000 0.00503005973506182\\
};
\end{axis}
\end{tikzpicture}%
\caption{Relative error}
\label{fig:relatverr}
\end{subfigure}
\caption{Comparison of full model and reduced order model for $mod-2476$.}
\label{fig:rom}
\end{figure}

In Figure.\ref{fig:rom}, the comparison of the full model and 100 dimensional ROM of $mod-2476$ is provided. From Figure.\ref{fig:sigmaplot} it is observed that the transfer functions of the full model and the ROM are identical. Figure.\ref{fig:abserror} and Figure.\ref{fig:relatverr} illustrate the absolute error and relative error of the ROM, respectively. From those figures, it is evident that the absolute error and relative error are significantly small.      

\subsection{Analysis of the eigenvalues}
\begin{figure}[]
\setlength{\figheight}{.4\linewidth}
\setlength{\figwidth}{.75\linewidth}
\centering
%
%
%
%
\begin{tikzpicture}

\begin{axis}[%
width=\figwidth,
height=\figheight,
scale only axis,
separate axis lines,
every outer x axis line/.append style={darkgray!60!black},
every x tick label/.append style={font=\color{darkgray!60!black}},
xmin=-0.04,
xmax=0.005,
xlabel={Real axis},
every outer y axis line/.append style={darkgray!60!black},
every y tick label/.append style={font=\color{darkgray!60!black}},
ymin=-0.107121663800708,
ymax=0.108391744340198,
ylabel={Imaginary axis},
axis on top,
legend columns=2,
legend style={nodes=right,anchor=south, at={(0.5,1.05)}}
]
\addplot [
color=blue,
only marks,
mark=o,
mark options={solid},
forget plot
]
table[row sep=crcr]{
0.002821871236311 0\\
-0.0346737242807555 0.0885105661130785\\
-0.0346737242807555 -0.0885105661130785\\
-0.0361793741298476 0.0901831195396975\\
-0.0361793741298476 -0.0901831195396975\\
-0.0014949263956072 0\\
-0.00285052483918429 0\\
-0.0190098755733593 0.0105430352491079\\
-0.0190098755733593 -0.0105430352491079\\
};
\addplot [
color=blue,
only marks,
mark=o,
mark options={solid},
forget plot
]
table[row sep=crcr]{
-0.00220778525702834 0\\
-0.0237846856093015 0\\
};
\addplot [
color=blue,
only marks,
mark=o,
mark options={solid},
forget plot
]
table[row sep=crcr]{
-0.0214518786671665 0\\
};
\addplot [
color=blue,
only marks,
mark=o,
mark options={solid},
forget plot
]
table[row sep=crcr]{
-1.01959676856992e-14 0\\
};
\addplot [
color=blue,
only marks,
mark=o,
mark options={solid},
forget plot
]
table[row sep=crcr]{
5.71112223980568e-14 0\\
6.48085695995449e-13 0\\
};
\addplot [
color=blue,
only marks,
mark=o,
mark options={solid},
forget plot
]
table[row sep=crcr]{
8.12324293178362e-13 0\\
};
\addplot [
color=blue,
only marks,
mark=o,
mark options={solid},
]
table[row sep=crcr]{
3.5515108720211e-12 0\\
};
\addlegendentry{original};
\addplot [
color=red,
only marks,
mark=+,
mark options={solid},
forget plot
]
table[row sep=crcr]{
0.00102298451444043 0\\
-0.00148528413190876 0\\
-0.0346737197845723 0.0885105739809589\\
-0.0346737197845723 -0.0885105739809589\\
-0.0361793831463839 0.0901831551647305\\
-0.0361793831463839 -0.0901831551647305\\
-0.00285052463278302 0\\
-0.0181636650289949 0.00753774727716888\\
-0.0181636650289949 -0.00753774727716888\\
-0.00220778525280713 0\\
2.52874294563364e-14 0\\
};
\addplot [
color=red,
only marks,
mark=+,
mark options={solid},
forget plot
]
table[row sep=crcr]{
-0.0237846856438811 0\\
1.36473187563237e-11 0\\
-0.0214518786323054 0\\
};
\addplot [
color=red,
only marks,
mark=+,
mark options={solid},
forget plot
]
table[row sep=crcr]{
3.95836313308765e-13 0\\
-8.59428165373449e-12 0\\
};
\addplot [
color=red,
only marks,
mark=+,
mark options={solid},
]
table[row sep=crcr]{
2.0303690790604e-12 0\\
};
\addlegendentry{stabilized};
\end{axis}
\end{tikzpicture}%
\caption{Eigenvalue stabilization of the model $mod-2476$.}
\label{fig:eigs}
\end{figure}
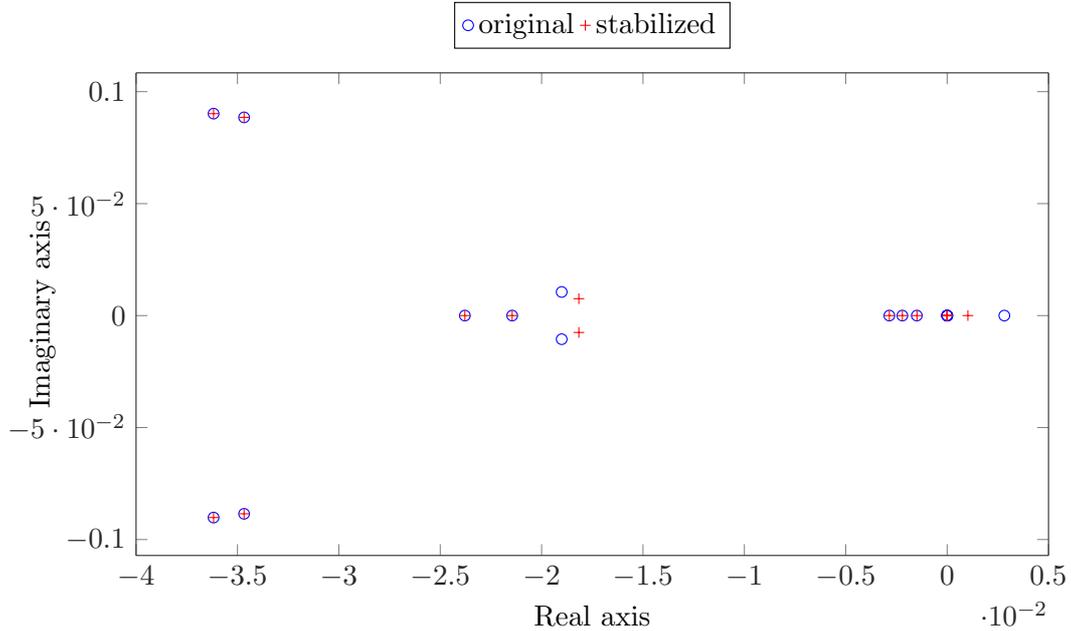

Figure.\ref{fig:eigs} depicts the magnified eigenvalues of the original system and the stabilized system. From the figure, it can be said that using the proposed technique the unstable eigenvalues of the original system can be sufficiently stabilized.  

\subsection{Stabilization of step-responses}
\begin{figure}[]
\setlength{\figheight}{.325\linewidth}
\begin{subfigure}{\linewidth}
\centering
\setlength{\figwidth}{.75\linewidth}
\input{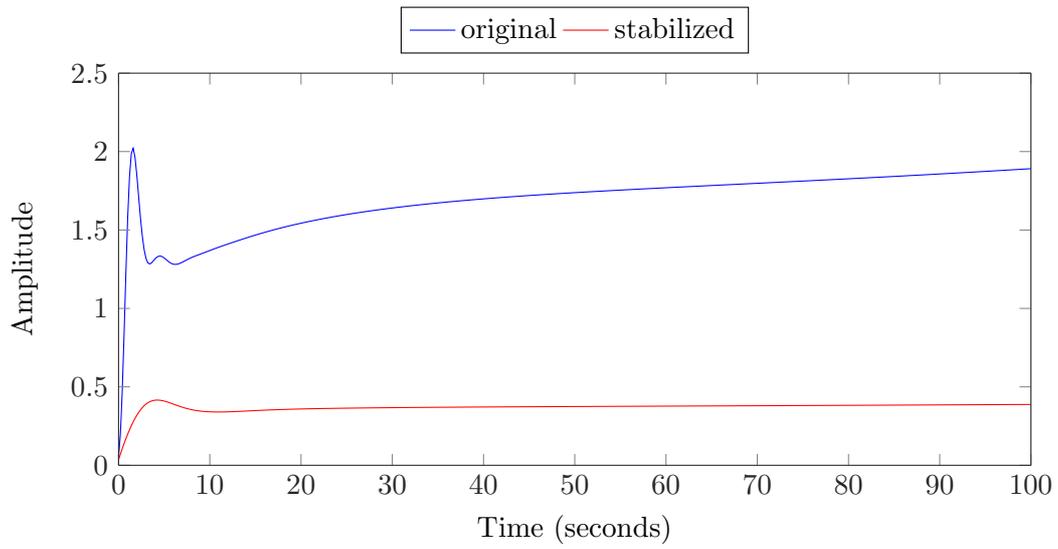}
\caption{Second input/first output}
\label{fig:step21}
\end{subfigure}
\begin{subfigure}{\linewidth}
\centering
\setlength{\figwidth}{.75\linewidth}
\input{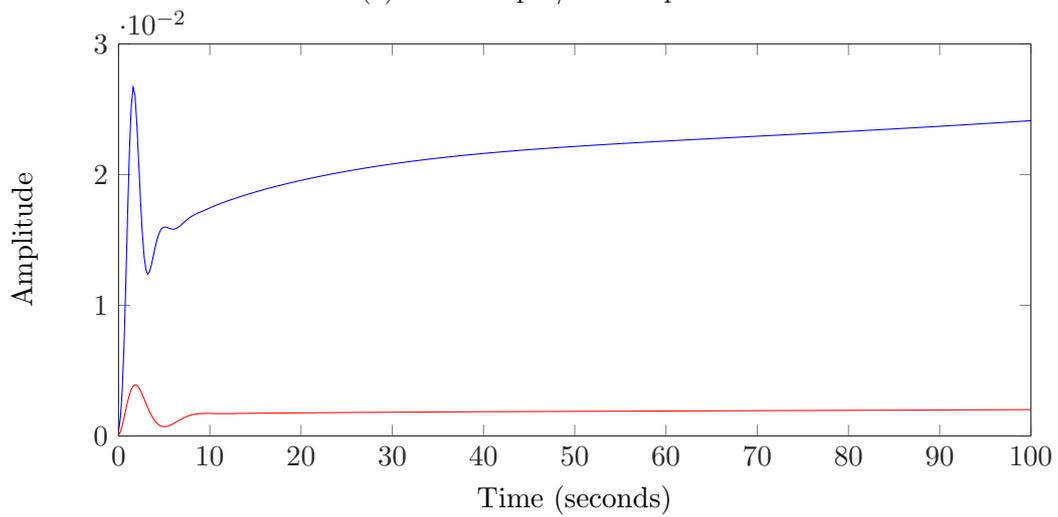}
\caption{Third input/second output}
\label{fig:step32}
\end{subfigure}
\begin{subfigure}{\linewidth}
\centering
\setlength{\figwidth}{.75\linewidth}
\input{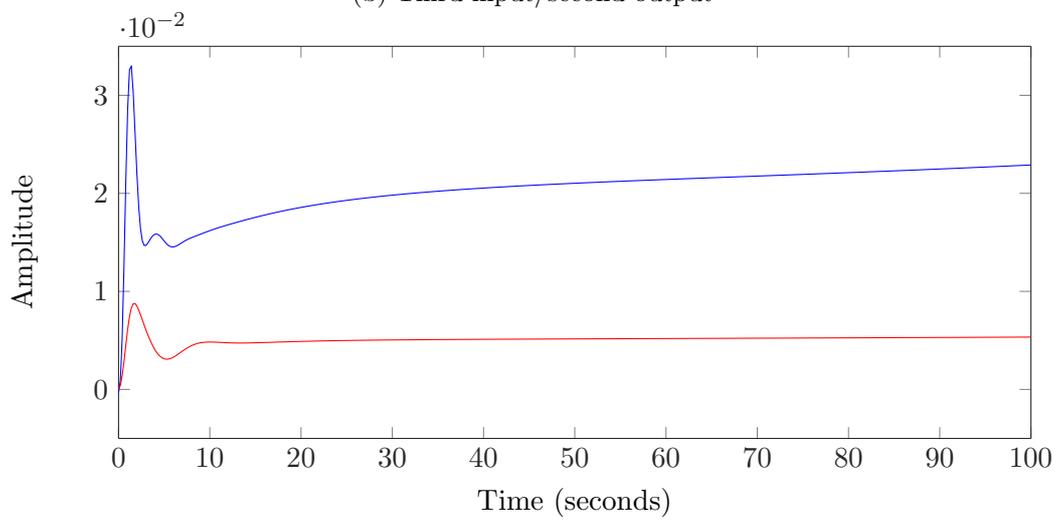}
\caption{Fourth input/third output}
\label{fig:step43}
\end{subfigure}
\caption{Stabilization of step-responses $mod-2476$.}
\label{fig:step_response}
\end{figure}

Figure.\ref{fig:step_response} represents the stabilization of the dominant step-responses of the target model. To reduce the size of the paper, only two step-responses are accounted. In Figure.\ref{fig:step21}, Figure.\ref{fig:step32}, and Figure.\ref{fig:step43} the step-responses of the original and stabilized systems are shown for second input/first output, third input/second output, and fourth input/third output, respectively. From the figurative comparison, it is revealed that the modified form of IRKA can be efficiently applied for the feedback stabilization of the unstable index-1 descriptor systems.

\subsection{Comparison with the previous works}
From the numerical computation, it has been seen that BIPS models described in Table.\ref{tab:BIPS_systems} show similarities of the stabilization of eigenvalues and step-responses for RKSM, LRCF-ADI based Newton-Kleinman (KN-LRCF-ADI), and IRKA techniques. On the other hand, the computation times differ significantly in those methods. In Table.\ref{tab:Comparison} the comparison of the simulation time of the mentioned three techniques for target models is displayed. 
\begin{table}[]
\centering
\caption{Time comparison of IRKA with RKSM and KN-LRCF-ADI techniques}
\label{tab:Comparison}
\begin{tabular}{cccc}\hline
Method &  RKSM   &   KN-LRCF-ADI    &   IRKA \\ \hline
$mod-606$ &  $1.81 \times 10^2$   &   $7.21 \times 10^2$    &   $4.11 \times 10^2$ \\ 
$mod-1998$ &  $1.41 \times 10^3$   &   $5.83 \times 10^3$    &   $4.89 \times 10^3$ \\ 
$mod-2476$ &  $3.06 \times 10^3$   &   $5.74 \times 10^3$    &   $4.92 \times 10^3$ \\ \hline
\end{tabular}
\end{table}

From Table.\ref{tab:Comparison}, it can be said that the IRKA approach is better than KN-LRCF-ADI in case of time comparison, whereas RKSM is the best.

\section{Conclusions}
In this work, a modified form of bi-tangential interpolation-based technique Iterative Rational Krylov Algorithm (IRKA) for the Riccati-based feedback stabilization of unstable index-1 descriptor systems, has been proposed and derived. The proposed method has been implemented to unstable power system models derived from Brazilian Interconnected Power System (BIPS). To pursue the desired goal, the Reduced-Order Model (ROM) based low-rank CARE has been solved by the MATLAB library command \texttt{care} and the corresponding low-rank feedback matrix has been estimated. The optimal feedback matrix for stabilizing the full model has been retrieved. The sparse form of the construction of the projection matrices and an approach to imply the initial feedback matrix for the treatment of the unstable systems have been discussed. As the validation process, the tabular and graphical approaches have been used for comparative analysis of the numerical results gained by the proposed technique. It has been observed that the ROMs attained by the proposed method properly represent the full model. From the stabilization of the eigenvalues and step-responses, it is evident that the proposed technique is efficient for the stabilization of unstable descriptor systems. Moreover, the computation time of the proposed technique with two previous techniques has been compared and found that it worked better than one of them but not the other. So, further research is needed to reduce the computation time of the proposed technique.   

\ack{This work is under the project \textit{"Computation of Optimal Control for Differential-Algebraic Equations (DAE) with Engineering Applications"}. This project is funded by United International University, Dhaka, Bangladesh. It starts from October 01, 2019, and the reference is IAR/01/19/SE/18.}

\section*{References}

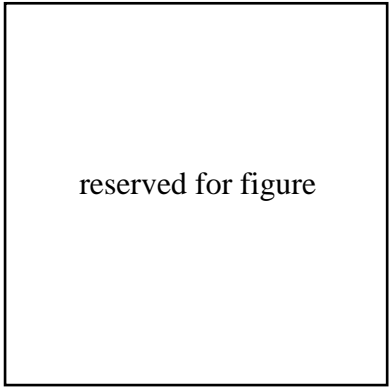
\end{document}